\documentclass[a4paper, 10pt, conference]{ieeeconf}
\IEEEoverridecommandlockouts
\usepackage[utf8]{inputenc} 
\usepackage{cite}

\usepackage{amsmath,amssymb,amsfonts}
\usepackage{graphicx}
\usepackage{textcomp}
\usepackage{xcolor}

\usepackage[utf8]{inputenc} 
\usepackage[T1]{fontenc}    
\usepackage{amsfonts}       
\usepackage{nicefrac}       

\usepackage{xcolor}         
\usepackage{mdframed}
\usepackage{tcolorbox}
\usepackage{amsmath,amssymb}

\usepackage{amsthm}



\usepackage[export]{adjustbox}
\usepackage{mathtools}
\usepackage{derivative}
\usepackage{array}
\usepackage{makecell}
\usepackage{longtable}
\usepackage{hhline}
\usepackage{nicematrix}
\usepackage[bottom]{footmisc}
\usepackage{multirow}
\usepackage{multicol}
\usepackage{booktabs}
\usepackage[mathscr]{euscript}
\usepackage{mathrsfs}
\usepackage{calligra}

\usepackage{tikz}
\usetikzlibrary{shapes,calc,arrows.meta,backgrounds,positioning,fit} 

\usepackage{parskip}
\usepackage{caption}
\usepackage[caption=false]{subfig}
\let\labelindent\relax
\usepackage[inline]{enumitem}

\usepackage{chngcntr}
\usepackage{empheq}

\usepackage[
  colorlinks=true,
  linkcolor=blue,
  citecolor=blue,
]{hyperref}
\usepackage{url}

\usepackage{bookmark}
\usepackage[capitalise,nameinlink]{cleveref}
\crefname{assumption}{Assumption}{Assumptions}
\Crefname{assumption}{Assumption}{Assumptions}

\usepackage{algorithm}
\usepackage{algorithmicx}
\usepackage{algpseudocode}
\algnewcommand\algorithmicoutput{\textbf{Output:}}
\algnewcommand\Output{\item[\algorithmicoutput]}

\newtheoremstyle{myplain}
{4pt}  
{4pt}  
{\itshape}  
{}      
{\bfseries} 
{.}     
{0.5em} 
{}
\theoremstyle{myplain}
\newtheorem{theorem}{Theorem}

\theoremstyle{definition}

\theoremstyle{remark}

\newcommand{\tr}[1]{\ensuremath{\mathrm{tr}\left[#1\right]}}

\newcommand{\E}[2][]{\ensuremath{\mathbb{E}_{#1}\left[ #2 \right]}}

\newcommand{\norm}[1]{\ensuremath{\left\| #1 \right\|}}

\DeclareMathOperator{\covMatrix}{\text{$\Sigma$}}

\def\BibTeX{{\rm B\kern-.05em{\sc i\kern-.025em b}\kern-.08em
    T\kern-.1667em\lower.7ex\hbox{E}\kern-.125emX}}

\begin{document}

\title{\LARGE \bf LQG solution for POMDP without estimating states: \\ A minimum variance approach}

\author{Ranjan Sarkar, Prabhat K. Mishra
	\thanks{Authors are with the Department of Artificial Intelligence at \break Indian Institute of Technology Kharagpur, West Bengal, India \break
    (\href{mailto:ranjan.sarkar.24@kgpian.iitkgp.ac.in}{\texttt{ranjan.sarkar.24@kgpian.iitkgp.ac.in}}, \href{mailto:pkmishra@ai.iitkgp.ac.in}{\texttt{pkmishra@ai.iitkgp.ac.in}}).}
}

\maketitle

\begin{abstract}
    This paper investigates the control of discrete-time linear time-invariant (LTI) systems subject to incomplete and corrupted measurements. Specifically, we focus on designing a Linear Quadratic Gaussian (LQG) controller without relying on explicit state estimation. By leveraging minimum variance duality, our approach allows the current control input to be represented as a linear function of available measurements and previously applied inputs, successfully reducing the task to a tractable deterministic optimization problem. We provide theoretical justification for this framework and demonstrate its practical effectiveness through numerical experiments.
\end{abstract}

\section{Introduction} \label{sec:intro}

The development of stochastic control theory has historically been anchored by the Separation Principle (or its special version, Certainty Equivalence) \cite{2013_separation, 2022_dual_Serdar}. Under this classical paradigm, the synthesis of an output feedback controller is strictly sequential: an optimal estimator (such as a Kalman filter) is designed first to reconstruct the hidden system states from noisy measurements, followed by the computation of a Linear Quadratic Regulator (LQR) or a feedback policy. The constrained stochastic optimal control methodologies \cite{hokayem2012stochastic, mishra2017output} also mandate the same two-stage sequential process. While this ``estimation-first'' approach provides a highly modular framework for idealized linear-Gaussian systems, modern structural analyses and engineering constraints reveal limitations in its universal application and ignited new interest in their study \cite{2020_LQR_Recht, tang2023analysis, 2026_LQG_games}. 

The assumption that estimation and control can be perfectly decoupled breaks down when systems deviate from unconstrained, continuous linear models \cite{2022_separation_distributed} except in a few special cases \cite{1999_separation_Khalil, 2025_separation_class, 2026_separation_Mesbah}.
The integration of digital communication channels inherently shatters certainty equivalence \cite{2023_separation_TCP}. As demonstrated by \cite{fu2012lack}, when a feedback channel is subject to quantization with a fixed bit rate, state estimation and control can no longer be fully separated. These limitations are most striking in systems with bilinear observation models, where the control input actively alters how the state is measured. In such settings, the Separation Principle fails entirely, and naively applying a standard, separated LQG controller can induce a loss of observability, and in certain cases inadvertently maximizing the quadratic cost rather than minimizing it \cite{sattar2025suboptimality}. 

To circumvent the limitations of explicit state estimation, an alternative paradigm \cite{bakolas2018constrained} focuses on directly parameterizing the control policy to be causal linear combinations of past output measurements and past inputs, the underlying stochastic optimal control problem can be elegantly transformed into a tractable, deterministic surrogate optimization problem. Given the structural sub-optimality and potential instability of the standard ``Estimation-first'' approach under constraints, there is a distinct need for an alternate methodology. Motivated by the results in minimum variance duality \cite{2022_mishra_automatica, kim2024variance, kim2024arrow}, this paper targets an alternate synthesis paradigm for linear systems subjected to Gaussian noise and corrupt measurements. Instead of estimating states as an intermediate step, our approach focuses on computing the LQR gain first. By mapping the statistical properties of the noise and the historical observation data directly into the cost function, we can embed this LQR gain into a deterministic optimization problem that directly yields the control action. 

The main contributions of this work are as follows:
\begin{itemize}
    \item A direct LQG control estimation approach for partially observable Markov decision process (POMDP) is introduced that eliminates the need for explicit state estimation.
    \item We utilize the minimum-variance duality \cite[Ch. 7.5]{astrom} with a suitable modification \cite{2022_mishra_automatica} required to incorporate the current measurement and derive the proposed method, established in \cref{thm:duality_cov}, by adopting a different initialization of the dual state from that used in \cite{2022_mishra_automatica}.
\end{itemize}

Our notations are standard and mostly defined before their first use. We let $\mathbb{R}$ denote the set of real numbers and $\mathbb{Z}$ the set of integers. The notations $\mathbb{Z}_{\geq 0}$ denotes the set of non-negative integers or natural numbers including $0$. A symmetric positive semi-definite (definite) matrix $M$ is denoted by $M \succeq 0$ ($M \succ 0$). The $n\times n$ identity matrix is denoted by $I_n$ or simply $I$ when dimensions are clear from the context. Similarly, $0_n$ denotes the zero matrix with all entries zero. A blockdiagonal matrix with diagonal entries $A$ and $B$ is defined by $\operatorname{blockdiag}(A,B)$.

We describe the problem setup in \S\ref{sec:background}, present our main result in \S\ref{sec:methodology}. All the experimental results are provided in \S\ref{sec:experiments} and we conclude in \S\ref{sec:conclusion}. Proofs of our main results are provided in the appendix.

\begin{figure*}[t]
    \centering
    
    \definecolor{hiddenfill}{RGB}{224,217,213}
    \definecolor{panelgray}{RGB}{237,237,237}
    \definecolor{darkteal}{RGB}{16,104,122}
    \definecolor{historytext}{RGB}{11,69,81}
    \definecolor{lightblue}{RGB}{202,238,251}
    \definecolor{darkline}{RGB}{64,64,64}
    \definecolor{dashgray}{RGB}{88,88,88}
    \definecolor{magentaDash}{RGB}{160,40,144}
    \definecolor{SlateBlue}{HTML}{2E4057}
    \definecolor{TealFill}{HTML}{A3CEF1}
    \definecolor{SumFill}{HTML}{E7ECEF}
    \definecolor{StateFill}{HTML}{E0D9D5}
    \definecolor{ObsFillLeft}{HTML}{4A5E6D}
    \definecolor{ObsFillRight}{HTML}{3A606E}
    \definecolor{HistoryOuter}{HTML}{9E2A2B}
    \definecolor{HistoryBg}{HTML}{F4F4F6}
    \definecolor{InnerBoxBg}{HTML}{E2EAF4}
    \definecolor{InnerBoxLine}{HTML}{2A6F97}
    
    \tikzset{
        >=Stealth,
        draw=SlateBlue,
        text=black,
        node distance=1cm and 1.2cm,
        circleNode1/.style={draw, circle, fill=SumFill, minimum size=0.9cm, outer sep=0.5pt, font=\large},
        circleNode2/.style={draw, circle, fill=StateFill, minimum size=1.05cm, outer sep=0.5pt, font=\large},
        dashedCircleNode/.style={draw, circle, densely dashed, fill=StateFill, minimum size=1.05cm, outer sep=0.5pt, font=\large},
        sumNode/.style={draw, circle, fill=SumFill, minimum size=1.05cm, font=\large, outer sep=0.5pt, text=SlateBlue},
        rectNode/.style={draw=none, rounded corners=3pt, fill=darkteal, text=white, minimum width=0.8cm, minimum height=0.8cm, inner sep=8pt, outer sep=0.5pt, font=\large},
        dashedrectNode/.style={draw, rounded corners=3pt, fill=TealFill, minimum width=1.2cm, minimum height=1cm, outer sep=0.5pt, inner sep=0.5pt, font=\large},
        controlCircleNode/.style={draw=InnerBoxLine, circle, fill=TealFill, minimum size=1.05cm, outer sep=0.5pt, font=\large},
        estimate/.style={circle, minimum size=1cm, inner sep=0pt, fill=TealFill, draw=darkline, outer sep=0.5pt, line width=1pt},
    }
    
    \subfloat[Standard\label{fig:standard}]{
    \begin{minipage}[c]{0.36\textwidth}
    \centering
    \resizebox{\linewidth}{!}{
    \begin{tikzpicture}[scale=1.0, transform shape]
            
        \node[circleNode2] (xt) {$\hat{x}_t$};
        \node[sumNode, right=of xt] (ut) {$\hat{u}_t$};
        \node[dashedrectNode, right=1.7cm of ut, inner sep=10pt] (sys) {System};
        \node[rectNode, below=of xt] (yt) {$y_0,\dots,y_t$};
        \node[rectNode, below=of sys] (yt1) {$y_{t+1}$};
        
        \draw[->, thick, SlateBlue] (xt) -- node[above=2pt, text=SlateBlue!80!black] {$K$} (ut);
        \draw[->, thick, SlateBlue] (ut) -- node[above=2pt, text=SlateBlue!80!black] {Apply} (sys);
        \draw[->, thick, SlateBlue] (sys) -- node[right=4pt, text=SlateBlue!80!black] {Measure} (yt1);
        \draw[->, thick, SlateBlue] (yt) -- node[right=2pt, text=black] {Estimate} (xt);
        \draw[->, thick, SlateBlue] (yt1) -- node[above=1pt]{Append} (yt);
    \end{tikzpicture}
    }
    \end{minipage}}
    \hspace{1.3cm}
    \subfloat[Proposed\label{fig:proposed}]{
    \begin{minipage}[c]{0.36\textwidth}
    \centering
    \resizebox{\linewidth}{!}{
    \begin{tikzpicture}[scale=1.0, transform shape]
            
    \node[sumNode] (ut) {$\hat{u}_t$};
    \node[dashedrectNode, right=2.2cm of ut, inner sep=10pt]
        (sys) {System};
    \node[rectNode, below=1.1cm of ut, xshift=0.0cm]
        (yt) {$y_{0:t}$};
    \node[rectNode, below=1.1cm of ut, xshift=-1.4cm]
        (ut_hist) {$u_{0:t-1}$};
    \node[rectNode, below=1.1cm of sys]
        (yt1) {$y_{t+1}$};

    \begin{scope}[on background layer]
        \node[
            draw=SlateBlue,
            rounded corners=4pt,
            fill=none,
            inner sep=8pt,
            fit=(ut_hist) (yt)
        ] (history) {};
    \end{scope}

    \draw[->, thick, SlateBlue]
        (ut) -- node[above=2pt, text=SlateBlue!80!black] {Apply} (sys);
    \draw[->, thick, SlateBlue]
        (sys) -- node[right=4pt, text=SlateBlue!80!black] {Measure} (yt1);
    \draw[->, thick, SlateBlue]
        (history.north) ++ (0.8,0) -- node[right=2pt, text=black] {Estimate} (ut);
    \draw[->, thick, SlateBlue]
        (yt1) -- node[above=1pt] {Append} (yt);

    \draw[->, thick, SlateBlue]
        (ut.west) -- ++(-0.856cm,0)
        -- node[pos=0.5, left=2pt] {Append}
        (ut_hist.north);
    \end{tikzpicture}
    }
    \end{minipage}}
    
    \caption{Standard approaches first estimate states then compute the control. Our proposed approach directly computes the control from the measurement data by utilizing minimum variance duality.}
    \label{fig:solution_architecture}
\end{figure*}

\section{Problem} \label{sec:background}

We consider a linear time-invariant (LTI) system with incomplete and corrupt measurements. The state $x_t \in \mathbb{R}^n$ is partially observable through lower-dimensional observations or measurements $y_t \in \mathbb{R}^m$, where $m < n$, affected by process noise $\xi_t$ and measurement noise $\omega_t$. We describe the measurement sequence with the control or input $u_t \in \mathbb{R}^{p}$ through following stochastic difference equations for $t \in \mathbb{Z}_{\ge 0}$,
\begin{subequations} \label{eq:system}
\begin{align}
    x_{t+1} &= A x_{t} + B u_{t} + \xi_{t}, \label{eq:state_evolution} \\
    y_{t} &= C x_{t} + \omega_{t}, \label{eq:output_evolution}
\end{align}
\end{subequations}
where
\begin{enumerate}[leftmargin = *, nosep, label=(1-\alph*), widest = b]
\item \label{asmp:initial state distribution}
    The initial state $x_0 \in \mathbb{R}^n$ is a random vector with mean $\hat{x}_0^-$ and covariance $\covMatrix^-_0 \succeq 0$.
\item \label{asmp:white noise}
    For all $t$, the noise $\xi_t$ and $\omega_t$ are identically distributed Gaussian random vectors with mean zero and covariances $Q \in \mathbb{R}^{n \times n}$, $R \in \mathbb{R}^{m \times m}$, respectively, where $Q \succeq 0, R \succ 0$ are symmetric positive semi-definite and positive definite matrices, respectively.
\item
    For each $t$, the initial state $x_0$, and the noise $\xi_t, \omega_t$ are mutually independent random vectors.
\item
    The pair $(A,C)$ is observable.
\end{enumerate}
We are interested in finding the control sequence $u_{0:t-1} := \{u_0, u_1, \ldots, u_{t-1}\}$, which minimizes the cost
\begin{equation}\label{eq:lqg_cost}
\mathrm{J}_{\text{LQG}} := \mathbb{E} \left[x_{t}^\top Q_t x_{t} + \sum_{i=0}^{t-1} \left(x_i^\top Q_c x_i + u_i^\top R_c u_i\right)\right],
\end{equation}
where expectation is with respect to the probability distribution of $\xi_0, \ldots, \xi_{t-1}$. The penalty matrices $Q_c, Q_t \succeq 0 \in \mathbb{R}^{n \times n}$ and $R_c \succ 0 \in \mathbb{R}^{p \times p}$ are given. The above problem \eqref{eq:lqg_cost} is known as \textit{Linear Quadratic Gaussian} (LQG) \cite{1995_Zervos_LQG}, and its closed form solution is a linear-feedback control $u_i = K_i x_i$ for $i=0,1,\dots,t-1$, when the states $x_i$ are perfectly measurable. The feedback gain
\begin{equation}\label{eq:K_gain_time}
    K_i = -(R_c + B^\top P_{i+1} B)^{-1} B^\top P_{i+1} A,
\end{equation}
is obtained with help of the backward Riccati recursion
\begin{align*}
    P_{i} = &~A^\top P_{i+1} A  + Q_c \\
        &- A^\top P_{i+1} B (R_c + B^\top P_{i+1} B)^{-1} B^\top P_{i+1} A,
\end{align*}
by initializing $P_{t} = Q_t$. We assume that
\begin{enumerate}[leftmargin = *, nosep, label = (A\arabic*)]
    \item \label{asmp:riccati_convergence} the pair $(A,B)$ is controllable and the pair $(A,D)$ is observable, where $D$ is such that $D^\top D = Q_c$.
\end{enumerate}
The Assumption \ref{asmp:riccati_convergence} allows us to compute $P_c \succeq 0$ by solving the discrete algebraic Riccati equation (DARE)
\begin{equation*}
    P_c = A^\top P_c A - A^\top P_c B (R_c + B^\top P_c B)^{-1} B^\top P_c A + Q_c,
\end{equation*}
due to \cite[p. 75]{bertsekas2005dynamic}.
Therefore, we can compute the feedback gain $K$ by
\begin{equation}\label{eq:K_gain}
    K = -(R_c + B^\top P_c B)^{-1} B^\top P_c A.
\end{equation}
Since we have access to the history of noisy and incomplete observations $y_{0:t}$ at time $t$ instead of $x_t$, we cannot compute $u_t = Kx_t$. The standard approach \cite[\S4.3]{bertsekas2005dynamic} first computes $\hat{x}_t$ via state estimator (e.g. Kalman filter) and then $\hat{u}_t = K\hat{x}_t$ (\cref{fig:standard}).
Instead, we aim to compute $\hat{u}_t$ directly from $y_{0:t}$ and $u_{0:t-1}$ (\cref{fig:proposed}).
The precise problem statement for this paper is as follows:
\begin{equation}
\left.
\begin{aligned}\label{eq:problem_setup}
    \underset{\hat{u}_t}{\text{minimize}} & \quad \mathbb{E}\left[ \lVert Kx_t - \hat{u}_t \rVert^2 \mid y_{0:t}, u_{0:t-1}\right] \\
    \text{subject to} & \quad  \text{POMDP } \eqref{eq:system},
\end{aligned}
~\right\}
\end{equation}
where $K$ is LQR gain \eqref{eq:K_gain}.

\section{Solution} \label{sec:methodology}

The minimum-variance duality \cite[Ch.~7.5]{astrom} provides $\hat{x}_t^{-}$ at time $t$ as a linear combination of measurements $y_0, \ldots, y_{t-1}$, by converting the state estimation problem into an equivalent optimal control problem with the help of a dual process. It is observed in \cite{2022_mishra_automatica} that the same approach can be used to include the most recent measurement $y_t$ to get $\hat{x}_t$ as a linear combination of all measurements $y_0, \ldots, y_t$ available at time $t$, by only initializing the dual process differently. In this article, we will initialize the dual process in such a way we can solve \eqref{eq:problem_setup}.  

Let us consider the dual process
\begin{subequations} \label{eq:dual_dynamics_1}
\begin{align}
  z_{i+1} &= A^\top z_i + C^\top \alpha_{i+1}, \quad i \in \mathbb{Z}_{\geq 0},
    \label{eq:dual_dynamics_post}\\
  z_0 &= K^\top + C^\top \alpha_0 . \label{eq:dual_post_init}
\end{align}
\end{subequations}

Defining the cost-per-stage $\ell(z_i,\alpha_i)=z_i^\top Q z_i +  \alpha_i^\top R \alpha_i$, we get the following result:
\begin{theorem} \label{thm:duality_cov}
The optimization problem \eqref{eq:problem_setup} is equivalent to
\begin{equation} \label{eq:opti}
\left.
\begin{aligned}
    \underset{\alpha_{0:t}}{\operatorname{min}} ~&~\tr{z_t^\top \covMatrix_0^- z_t + \alpha_t^\top R \alpha_t + \sum_{i=0}^{t-1} \ell(z_i,\alpha_i)} \\
    \text{\upshape s.t.}
    ~&~\text{\upshape dual dynamics}~\eqref{eq:dual_dynamics_1}.
\end{aligned}
\right\}
\end{equation}
The minimizer $\hat{u}_t = \E{u_t \mid y_{0:t}, u_{0:t-1}}$ is given by
    \begin{equation}\label{eq:estimated_control}
        \hat{u}_t = z_t^\top \hat{x}_0^- + \sum_{i=0}^{t-1} z_i^\top B {u}_{t-i-1} - \sum_{i=0}^{t} \alpha_{i}^\top y_{t-i}.
    \end{equation}
\end{theorem}

A proof of \cref{thm:duality_cov} is given in the Appendix and a pseudo-code is provided in \cref{alg:uncons}. We have the following interesting observations related to \cref{thm:duality_cov}.
\begin{enumerate}[leftmargin = *, nosep, label = R\arabic*.]
    \item The optimization problem \eqref{eq:problem_setup} is stochastic but its equivalent \eqref{eq:opti} is purely deterministic. 
    \item If we substitute $K=I$ in \eqref{eq:opti}, we can get state estimate $\hat{x}_t$ according to \cite[Lemma 1]{2022_mishra_automatica}. 
    \item If we substitute $K=I$ and $\alpha_0 = 0$, we get the results of \cite[Ch.~7.5]{astrom}.
    \item Theorem \ref{thm:duality_cov} demonstrates that the optimization problem \eqref{eq:opti} can be used to estimate any linear function of $x_t$ without estimating $x_t$ by modifying $z_0$ accordingly. 
    \item Our class of feedback policies \eqref{eq:estimated_control} is the same as in \cite{bakolas2018constrained}.
    \item \label{rem:constraint} We can get constrained LQG control from \eqref{eq:opti} by simply including the affine constraint $\hat{u}_t \in \mathcal{U}$ in \eqref{eq:opti}. For convex $\mathcal{U}$, the resulting optimization problem remains convex because the inverse image of a convex set under an affine function is convex \cite[\S2.3.2]{boyd2004convex}.
\end{enumerate}

\begin{algorithm}[htbp]
\caption{Computation of LQG Control}
\label{alg:uncons} 
\begin{algorithmic}[1]
\Require $A, B, C, Q_c, R_c, Q, R, \hat{x}^-_0, \covMatrix^-_0$.
\State Compute the feedback gain $K$ by \eqref{eq:K_gain}.
\For{\textbf{each} $t$}
    \State Measure $y_t$.
    \State Solve the optimization \eqref{eq:opti}.
    \State Compute $\hat{u}_t$ by \eqref{eq:estimated_control}.
    \State Set $u_t = \hat{u}_t$ and apply $u_t$ to the system \eqref{eq:system}.
\EndFor
\end{algorithmic}
\end{algorithm}

Problem~\eqref{eq:opti} is a finite-horizon LQR problem for the dual process \eqref{eq:dual_dynamics_post}, except that the initial condition~\eqref{eq:dual_post_init}
depends on the decision variable $\alpha_0$. 
A straight-forward application of dynamic programming gives us the following result:

\begin{theorem}\label{thm:dual_solution}
Let $S_t = \Sigma_0^-$ and for $i = t-1, \dots, 0$,
\begin{equation}\label{eq:dual_riccati}
    S_i = Q + A S_{i+1} \big( I - C^\top L_{i+1} \big) A^\top ,
\end{equation}
where $L_i := (R + C S_i C^\top)^{-1} C S_i$. Then the minimizer of \eqref{eq:opti} is given by
\begin{equation}\label{eq:dual_gains}
\begin{aligned}
    \alpha_0^{\star} &= -L_0 K^\top, \\
    \alpha_{i+1}^{\star} &= -L_{i+1} A^\top z_i^{\star}, \quad i = 0, \dots, t-1,
\end{aligned}
\end{equation}
where $z_{0:t}^{\star}$ is obtained by substituting $\alpha_{0:t}^{\star}$ in \eqref{eq:dual_dynamics_1}.
\end{theorem}
A proof of \cref{thm:dual_solution} is provided in the Appendix.

\section{Numerical Experiments} \label{sec:experiments}
In this section, we evaluate our proposed method on two systems: an unconstrained cart-pole, and a constrained quadruped robot. In both cases, we compare with the standard approach of \cref{fig:standard}, i.e., a Kalman filter followed by $\hat u_t = K\hat x_t$, and use the noise-free perfect state feedback $u_t = Kx_t$ as a reference. In particular, the perfect state feedback case refers to $\xi_t, w_t = 0$ and $C=I$. All simulations are performed in Python~3.12.3 using \textit{CasADi} \cite{andersson2018casadi} with the \textit{IPOPT} solver.

\subsection{Experiment 1: Unconstrained Setting} \label{exp1}

We consider the inverted pendulum on a cart, linearized about the upright equilibrium, with the model and physical parameters as in \cite{ctms_matlab}. The state $x_t = [p_t\ \ v_t\ \ \theta_t\ \ \Omega_t]^\top$ collects the cart position, velocity, the pendulum angle and angular velocity, the input $u_t \in \mathbb{R}$ is the horizontal force on the cart, and only $p_t$ and $\theta_t$ are measured, i.e., $y_t = [p_t\ \ \theta_t]^\top \in \mathbb{R}^2$. The model is discretized by zero-order hold with $\Delta t = 0.01$\,s. We set $Q = 0.01I_4$, $R = 0.1I_2$, $Q_c = 100\,I_4$ and $R_c = 0.1$, which give $K = [15.05\ \ 23.62\ \ {-113.23}\ \ {-26.97}]$ by \eqref{eq:K_gain}.
The system starts from an unstable state $x_0 \sim \mathcal{N}(\hat x_0^-, \Sigma_0^-)$, where $\hat x_0^- = [-1\ \ 0\ \ 0.15\ \ 0]^\top$ and $\Sigma_0^- = 2I_4$.

\begin{figure}[htbp]
    \centering
    \includegraphics[width=\linewidth]{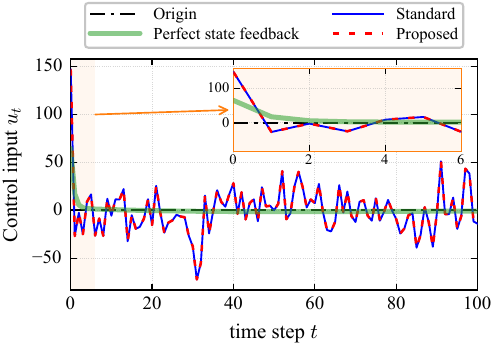}
    \caption{LQG control converges to zero when  states are perfectly measurable. Both the standard and the proposed approaches result in the same control sequence for the POMDP \eqref{eq:system}.}
    \label{fig:control}
\end{figure}

\begin{figure}[htbp]
    \centering
    \includegraphics[width=\linewidth]{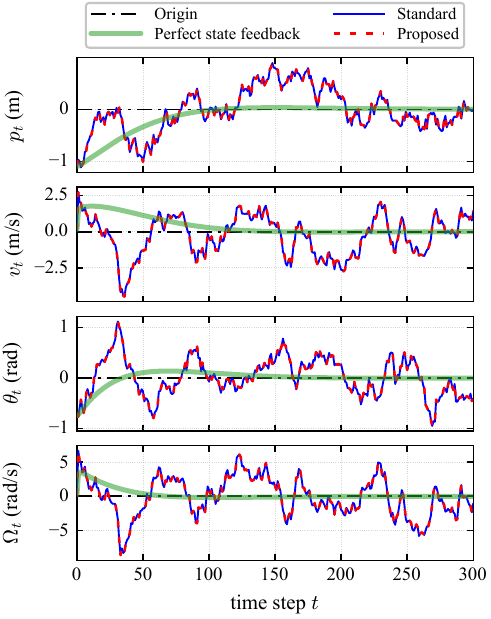}
    \caption{Both standard and proposed approaches give identical state trajectories, fluctuating around the noise-free state trajectory obtained by perfect state-feedback control.}
    \label{fig:state_dim}
\end{figure}

In \cref{fig:control}, we compare our proposed approach (red dashed) with standard approach (blue), and observe that both of them provide the same control sequence.
The corresponding state trajectories are presented in \cref{fig:state_dim}.

\subsection{Experiment 2: Constrained Setting
} \label{exp2}

\begin{figure}[H]
    \centering
    \includegraphics[width=0.9\linewidth]{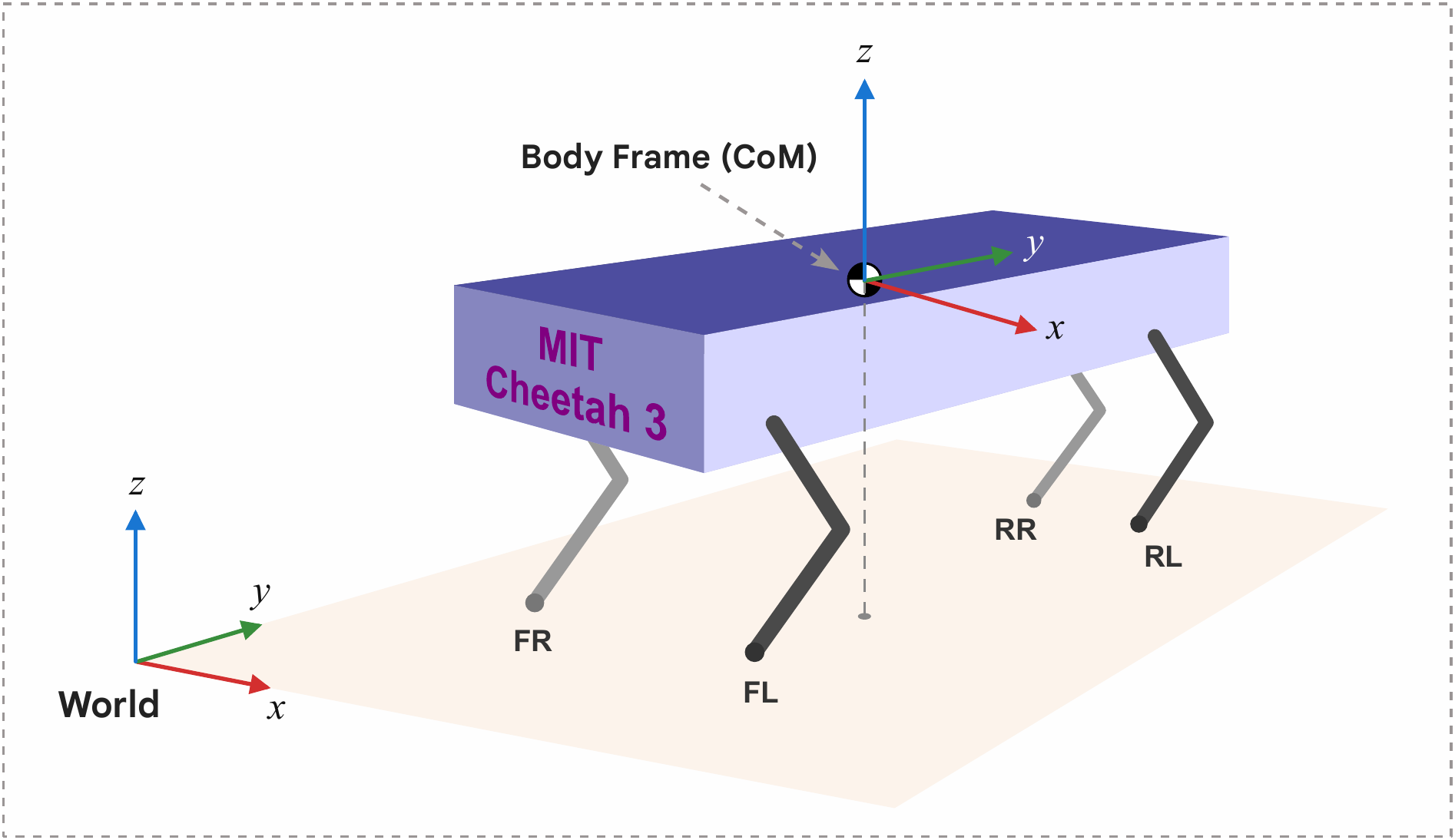}
    \caption{Quadruped model in the world coordinate system \cite{2018_mit_cheetah3}}
    \label{fig:quadruped}
\end{figure}

We consider the simplified rigid body dynamics of the quadruped robot MIT Cheetah 3 (\cref{fig:quadruped}) described in \cite{2018_mit_cheetah3} modeled as a discrete-time linear dynamical system.
The state vector $x_t \in \mathbb{R}^{13}$ at time-step $t$ completely characterizes the 3D translation and orientation dynamics of the system and is defined as:
\begin{equation*}
    x_t = \begin{bmatrix} \Theta_t^\top & p_t^\top & \omega_t^\top & v_t^\top & g \end{bmatrix}^\top \in \mathbb{R}^{13},
\end{equation*}
where $\Theta_t = [\phi_t\ \ \theta_t\ \ \psi_t]^\top \in \mathbb{R}^3$ represents the body orientation in Z-Y-X Euler angles (roll, pitch, yaw), $p_t \in \mathbb{R}^3$ is the position of the center of mass (CoM) in $(X,Y,Z)$ world coordinates, $\omega_t \in \mathbb{R}^3$ is the angular velocity, $v_t \in \mathbb{R}^3$ is the linear velocity, and $g = - 9.8$ is the constant gravity parameter.
The control input vector $u_t \in \mathbb{R}^{12}$ contains the 3 dimensional ground reaction forces applied by each of the four legs:
$$
u_t = \begin{bmatrix}
    f_{1,t}^\top & f_{2,t}^\top & f_{3,t}^\top & f_{4,t}^\top
\end{bmatrix}^\top \in \mathbb{R}^{12},
$$
where $f_{i,t} = [F_x\ \ F_y\ \ F_z]^\top \in \mathbb{R}^3$ represents the force components for the $i^\text{th}$ leg. The measurement vector $y_t \in \mathbb{R}^{6}$ contains noisy measurements of the orientation $\Theta_t$ and
the CoM position $p_t$. The continuous-time model of~\cite[Eq.~16]{2018_mit_cheetah3} is discretized exactly with $\Delta t = 0.03$\,s, the physical parameters are taken from \cite[Table~I]{2018_mit_cheetah3}.

The prior mean $\hat x_0^-$ and the covariance $\Sigma_0^-$ are listed in \cref{tab:cheetah},
and the true initial state is drawn as $x_0 \sim \mathcal{N}(\hat x_0^-, \Sigma_0^-)$. 
The noise covariances are given as
\begin{align*}
    Q &= \operatorname{blockdiag}((0.005)^2 I_6, (0.05)^2 I_3, (0.03)^2 I_3,0), \\
    R &= (0.05)^2 I_6.
\end{align*}
The LQR weights are $R_c = 0.01\,I_{12}$ and
\begin{equation*}
  Q_c = 50\operatorname{blockdiag}\big(10 I_3, \operatorname{diag}(10,10,50), I_3, I_3, 100\big),
\end{equation*}
and $K$ is computed from the DARE for the 12 controllable states. The task is to regulate the robot to the standing reference $x_{\mathrm{ref}} = [0_{1\times5}\ \ 0.4\ \ 0_{1\times6}\ \ {-9.8}]^\top$, i.e., a CoM height of $0.4$\,m, over 200 time steps. Each foot in MIT Cheetah dynamics is subject to the friction pyramid $|F_x| \le \mu F_z$, $|F_y| \le \mu F_z$ with $\mu = 0.5$ and to the vertical force limits $20 \le F_z \le 250$\,N, so $\mathcal{U}$ is a polytope. Therefore, we included constraints in \eqref{eq:opti} as noted in Remark \ref{rem:constraint}. For standard approach, a saturation function with desired bounds is applied to $K\hat x_t$ to satisfy constraints.

\begin{table}[t]
  \centering
  \caption{Initial state and covariances of the quadruped experiment (\S\ref{exp2})}
  \label{tab:cheetah}
  \setlength{\tabcolsep}{4pt}
  \begin{tabular}{l c c c c}
    \hline
    Block & $\hat x_0^-$ & $\Sigma_0^-$ \\
    \hline
    $\Theta$ (rad)      & $(0.25, -0.25, 0.30)$ & $(0.05)^2 I_3$ \\
    $p$ (m)             & $(0.40, -0.35, 0.28)$ & $(0.05)^2 I_3$ \\
    $\omega$ (rad/s)    & $(1.0, -1.0, 0.5)$    & $(0.1)^2 I_3$  \\
    $v$ (m/s)           & $(0.5, -0.5, -1.2)$   & $(0.1)^2 I_3$  \\
    $g$ (m/s$^2$)       & $-9.8$                & $0$ \\
    \hline
  \end{tabular}
\end{table}

\begin{figure}[htbp]
    \centering
    \includegraphics[width=\linewidth]{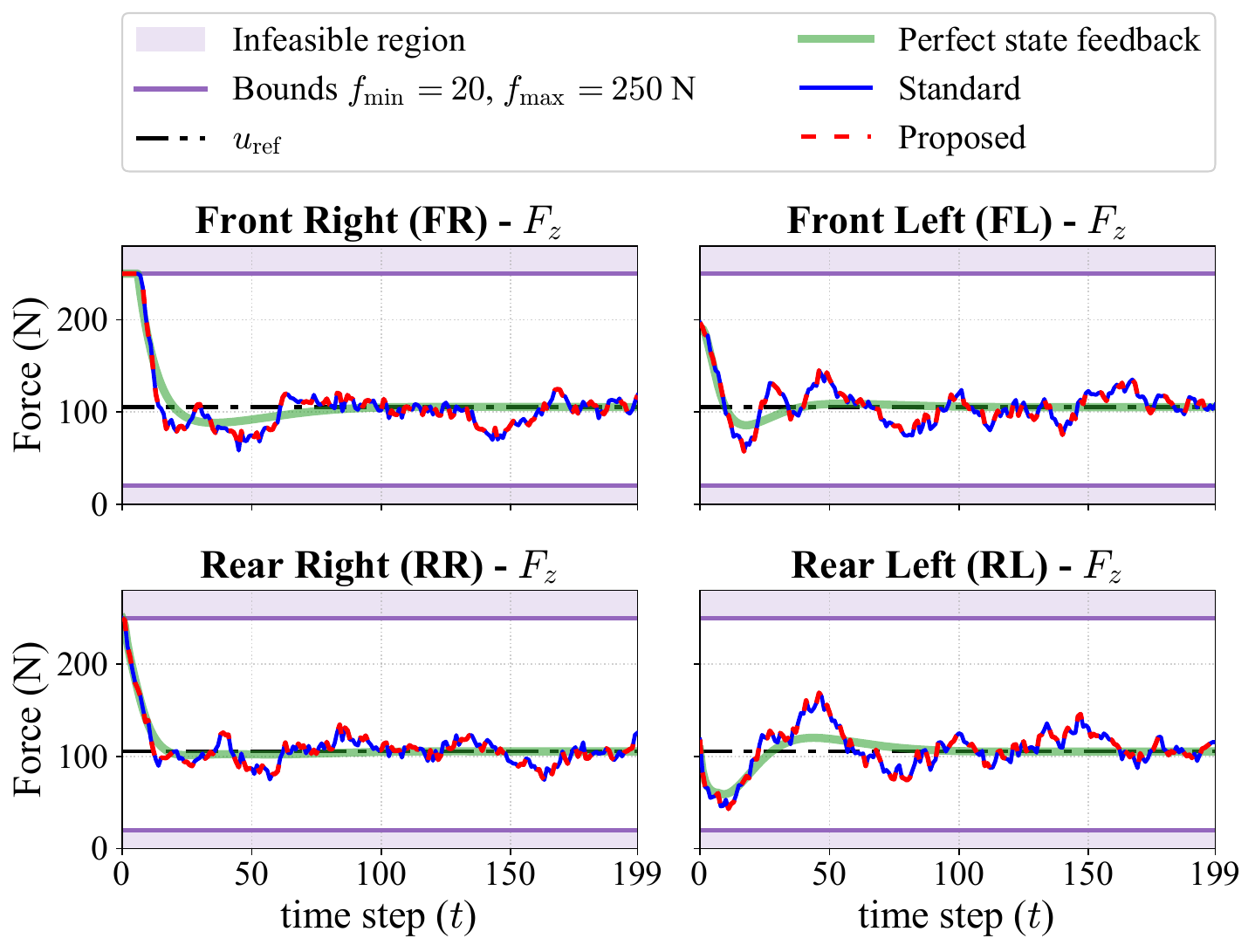}
    \caption{The perfect state-feedback control $F_z$ converge to their steady-state reference values $u_{\text{ref}}$, whereas the forces under both the standard and proposed methods approach these values but continue to oscillate due to persistent stochastic disturbances.}
    \label{fig:control2}
\end{figure}

\begin{figure}[htbp]
    \centering
    \includegraphics[width=\linewidth]{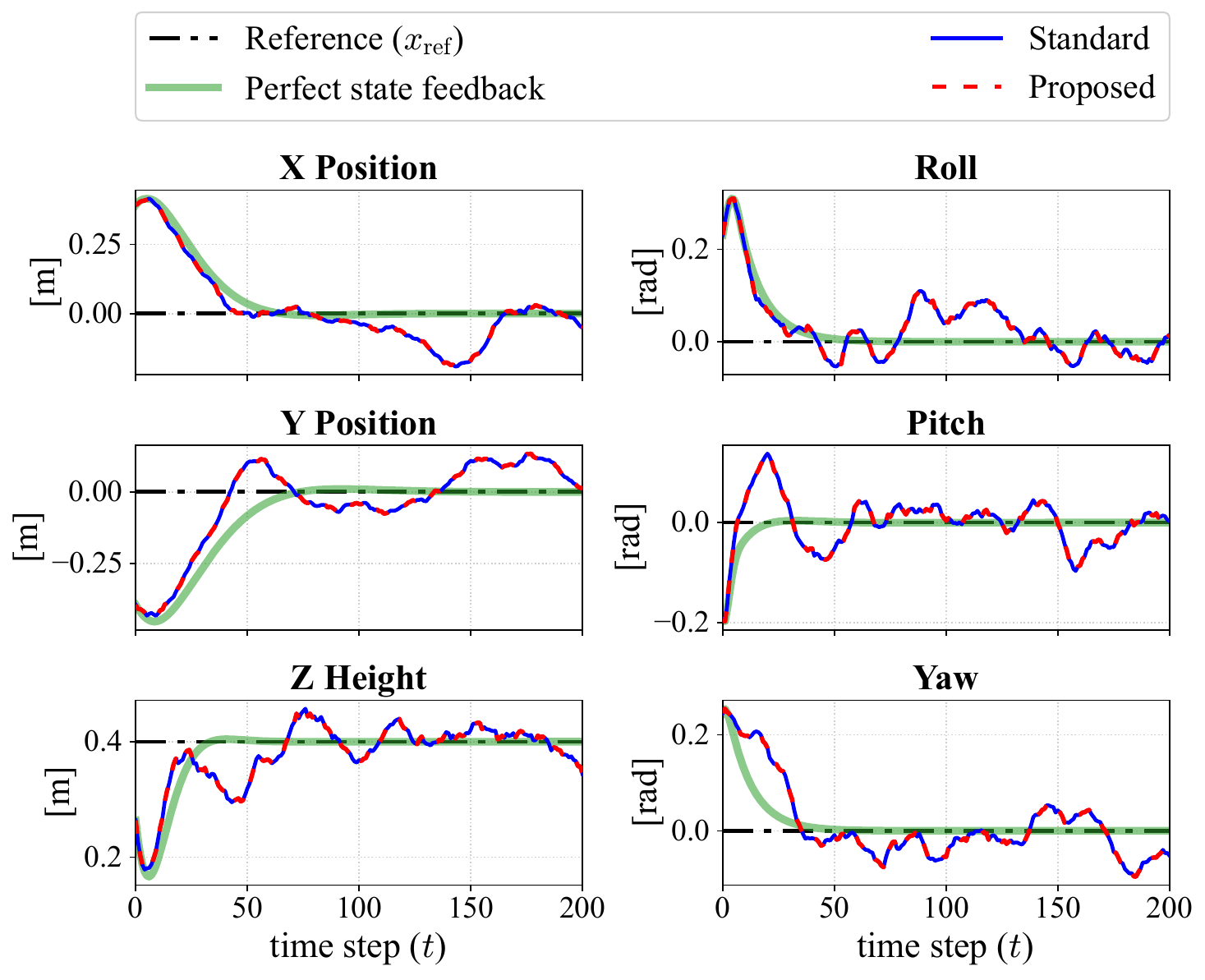}
    \caption{Both approaches exhibit similar state trajectories, oscillating around the trajectory obtained by the perfect state-feedback.}
    \label{fig:state2}
\end{figure}

\cref{fig:control2} shows the vertical forces and \cref{fig:state2} the CoM position and orientation. Both approaches keep the forces within their bounds and oscillate around the perfect state-feedback trajectories due to the noise. 


\section{Conclusion} \label{sec:conclusion}
We proposed an approach to directly estimate LQG control for POMDP without relying on explicit state estimation. By leveraging the minimum-variance duality, we developed an optimization-based approach that infers control inputs from noisy observations over a finite horizon. The proposed method offers a conceptually distinct alternative to the classical separation-based architecture and reduces dependence on intermediate state estimation.

\section*{Acknowledgment}
This work is supported in part by the ANRF PMECRG grant ANRF/ECRG/2024/004853/ENS and in part by the Faculty Startup Research Grant IIT/SRIC/AI/ODS/2024-2025/156. Authors are thankful to Prof. Prashant Mehta and Prof. Jin Won Kim for initial discussions related to minimum variance duality.

\appendix \label{appendix}

\begin{proof}[Proof of \cref{thm:duality_cov}]
    Using the relation $z_0 - C^\top \alpha_0 = K^\top$, we can write
    \begin{align}
        u_t = K x_t
        &= (z_0 - C^\top \alpha_0)^\top x_t \notag \\
        &= z_0^\top x_t - \alpha_0^\top C x_t. \label{eq:u_t}
    \end{align}
    Now using the system dynamics \eqref{eq:system} and dual dynamics \eqref{eq:dual_dynamics_post}, we obtain
    \begin{align}
        z_i^\top x_{t-i} &= z_i^\top (A x_{t-i-1} + B u_{t-i-1} + \xi_{t-i-1}), \label{eq:term1} \\
        z_{i+1}^\top x_{t-i-1} &= (A^\top z_i + C^\top \alpha_{i+1})^\top x_{t-i-1}. \label{eq:term2}
    \end{align}
    Subtracting \eqref{eq:term2} from \eqref{eq:term1} and summing over $i=0,\dots,t-1$ yields the telescoping sum as
    \begin{align*}
        z_0^\top x_t &= z_t^\top x_0 + \sum_{i=0}^{t-1} \left( z_i^\top x_{t-i} - z_{i+1}^\top x_{t-i-1} \right) \\
        &= z_t^\top x_0 + \sum_{i=0}^{t-1} \left( z_i^\top ( B u_{t-i-1} + \xi_{t-i-1}) - \alpha_{i+1}^\top C x_{t-i-1}\right).
    \end{align*}
    Substituting the above expression of $z_0^\top x_t$ in \eqref{eq:u_t}, we get
    \begin{equation*}
        u_t = z_t^\top x_0 + \sum_{i=0}^{t-1} z_i^\top (B u_{t-i-1} + \xi_{t-i-1}) -\sum_{i=0}^{t} \alpha_{i}^\top C x_{t-i}.
    \end{equation*}
    Since $\hat{u}_t := \E{u_t \mid {y}_{0:t}, u_{0:t-1}}$, we get \eqref{eq:estimated_control}.
    Subtracting \eqref{eq:estimated_control} from the above expression, we obtain the estimation error
    \begin{equation*}
        (u_t - \hat{u}_t) = z_t^\top (x_0 - \hat{x}_0^-) + \sum_{i=0}^{t-1} z_i^\top \xi_{t-i-1} + \sum_{i=0}^{t} \alpha_i^\top \omega_{t-i}.
    \end{equation*}
    Therefore, given $y_{0:t}, u_{0:t-1}$, the cost in \eqref{eq:problem_setup}
    \begin{equation*}
        \E{\norm{u_t - \hat{u}_t}^2} = \tr{\E{(u_t-\hat{u}_t)(u_t-\hat{u}_t)^\top}}
    \end{equation*}
    is simplified to the expression
    \begin{align}
        \tr{z_t^\top \covMatrix_0^- z_t + \sum_{i=0}^{t-1} z_i^\top Q z_i + \sum_{i=0}^{t} \alpha_i^\top R \alpha_i},
    \end{align}
where cross terms vanish due to independence and zero-mean assumptions of the noise terms \ref{asmp:white noise}.
\end{proof}

\begin{proof}[Proof of \cref{thm:dual_solution}]
    We use dynamic programming to solve the optimization problem \eqref{eq:opti}.
    For a given dual state $z_i=z$, $\alpha_i$ is independent of $z$. Therefore, we define the optimal cost-to-go as
    \[
    V_i(z) :=
    \min_{\alpha_{i+1:t}} ~\operatorname{tr}\!\left[
    \sum_{j=i}^{t-1} \left(z_j^\top Qz_j+\alpha_{j+1}^\top R\alpha_{j+1}\right) + z_t^\top\Sigma_0^-z_t
    \right],
    \]
    subject to the dual dynamics \eqref{eq:dual_dynamics_1}. At the terminal time,
    $V_t(z)=\operatorname{tr}(z^\top\Sigma_0^-z)$, so $S_t=\Sigma_0^- \succeq 0$. Suppose $V_{i+1}(z)=\operatorname{tr}(z^\top S_{i+1}z)$ with $S_{i+1} = S_{i+1}^\top \succeq 0$.
    The Bellman recursion gives
    \begin{equation} \label{eq:expanded_bellman}
        \begin{aligned}
        V_i(z) &= \min_{\alpha_{i+1}}~\operatorname{tr} \Bigl[
        z^\top Qz + \alpha_{i+1}^\top R\alpha_{i+1}\\
        &+(A^\top z+C^\top\alpha_{i+1})^\top
        S_{i+1}(A^\top z+C^\top\alpha_{i+1})
        \Bigr].
        \end{aligned}
    \end{equation}
    Differentiating the cost in \eqref{eq:expanded_bellman} with respect to $\alpha_{i+1}$ at $z=z_{i}$, using the properties of trace derivatives \cite[Prop.~10.7.4]{bernstein2009matrix} and setting the derivative to zero gives
    \begin{align}
        (R + CS_{i+1} C^\top)\alpha_{i+1} + C S_{i+1} A^{\top} z_i = 0.
    \end{align}
    Since, $R \succ 0$, so $(R+CS_{i+1}C^\top) \succ 0$ is invertible. Therefore, we have
    \begin{align*}
        \alpha_{i+1}^{\star} &= -L_{i+1}A^\top z_i, \\
        L_{i+1} &= (R+CS_{i+1}C^\top)^{-1}CS_{i+1}.
    \end{align*}
    Substituting the minimizer into the Bellman recursion gives $V_i(z)=\operatorname{tr}(z^\top S_i z)$, with
    \[
    S_i=Q+AS_{i+1}(I-C^\top L_{i+1})A^\top = S_i^\top.
    \]
    Thus, the backward recursion determines $S_{0:t}$ and the optimal controls $\alpha_{1:t}^{\star}$. To get optimizer $\alpha_0^{\star}$, we solve
    \[
    \min_{\alpha_0}~\tr{
    \alpha_0^\top R\alpha_0 + (K^\top+C^\top\alpha_0)^\top S_0(K^\top+C^\top\alpha_0)}.
    \]
    Setting its gradient to zero gives $$(R+CS_0C^\top)\alpha_0 + CS_0K^\top = 0,$$ and therefore
    \[
    \alpha_0^{\star}=-L_0K^\top.
    \]
\end{proof}

\bibliographystyle{ieeetr}
\bibliography{ref}

\end{document}